\documentclass [12pt]{article}
 \usepackage{graphicx,amssymb}
 \usepackage{epsfig}
\begin{document}

\centerline{\large \bf On a New Mechanism of Pattern Formation}

\vspace*{0.25cm}

\centerline{\large \bf in Population Dynamics}

\vspace*{1cm}

\centerline{S. Genieys$^1$, V. Volpert$^1$, P. Auger$^2$}

\vspace*{0.25cm}

\centerline{$^1$ Camille Jordan Institute of Mathematics, UMR 5208
CNRS, University Lyon 1} \centerline{69622 Villeurbanne, France}

\centerline{$^2$ Institute of Research and Development, 93143 Bondy, France}

\vspace*{2cm}

{\bf Abstract.}
We study a reaction-diffusion equation with an integral term describing 
nonlocal consumption of
resources. We show that a homogeneous equilibrium can lose its stability 
resulting in appearance
of stationary spatial structures. It is a new mechanism of pattern formation 
in population dynamics
that can explain emergence of biological species due to intra-specific 
competition and random mutations.
Travelling waves connecting an unstable homogeneous equilibrium and a 
periodic in space
stationary solution are studied numerically.

\vspace*{2cm}


 \setcounter{section}{1}
 \setcounter{equation}{0}

 \centerline{\bf 1. Introduction}

 \vspace*{0.5cm}

The concept of degeneracy in biology was proposed by Edelman in 1978
in relation with natural selection in biological systems (see [1]).
Degeneracy, according to his definition, means the ability of elements
that are structurally different to perform the same function or
yield the same output. He points out that degeneracy is necessary
for natural selection where natural selection is understood here
not only in the context of the evolution of biological species but
for many other biological systems including cell populations,
immune systems, brain cortex and others.
Degeneracy is related to robustness of biological systems and to
their ability to compensate lost or damaged functions.
For example, some proteins thought to be
indispensable for human organisms can be completely absent in
some individuals (see [1]).
There are also many examples related to the cortex activity where
some functions of the damaged area can be taken by some other areas.

In [2], Kupiec and Sonigo emphasize many physiological situations
such as cell differentiation or immunology, where darwinian selection
processes arise. They show that degeneracy and competition are very
common phenomena leading to the emergence of structures. 

Atamas studies the degeneracy in terms of signals and recognizers
[3]. Each recognizer can receive not only the signal, which
corresponds it exactly, but also close signals. He suggests that
each element recognizes a corresponding signal with a very high
specificity, which is not however absolute. The same element can
recognize other more or less similar signals though with less
probability. For example, it can be antigen receptors and
antibodies that have relative but not absolute specificity of
recognition. He fulfils numerical simulations with the individual
based approach where the recognizers can receive signals with
degenerate specificity and can reproduce themselves. This can be
considered as a model for natural selection with degenerate
recognition. The results of the numerical simulations show that
initially localized population of recognizers can split into two
subpopulations. Therefore the degeneracy can result in the
emergence of new species due to intraspecific competition and
selection.


 In this work we study degeneracy in biological system with
 continuous models.
 We consider them in the context of population dynamics.
 Starting with the classical logistic equation describing
 reproduction and migration of individuals of a biological species,
 we introduce a nonlocal consumption of resources that can
 influence the death function.
 One of possible interpretation of this model is as follows.
 Suppose that each individual is characterized by his positions
 $(x,y)$ at time $t$.
 The characteristic time scale is related to reproduction
 and can be measured for example in years.
 The individual can slowly migrate changing from time to time
 his position.
 For example it can be related to a new nesting place every
 year or every several years.
 The reproduction takes place at the same space point.
 However, contrary to conventional models in population dynamics
 we consider the space point $(x,y)$ not as the exact position
 of the individual but as his average position.
 This means that his nesting (reproduction) place is still a
 point but he can consume resources in some area around the
 point $(x,y)$.
 Therefore the death function, which takes into account the
 competition for the resources, depends not only on the
 individuals located at the point $(x,y)$ but also on those who
 are in some area around.
 We can ask the question, why the individuals coming to
 consume the resources at the point $(x,y)$ influence the death
 function at this point.
 One of possible answer could be for example that young
 individuals who are attached to their nesting place and
 who cannot yet consume the resources in a bigger area will have
 less chances to survive.
 We will specify the model in the next section.

 We will discuss also another interpretation of the model related
 to the evolution of biological species.
 In this case the space variable corresponds to some morphological
 characteristics and diffusion to a random small change of this
 characteristics due to mutations.
 The degeneracy in this case means that individuals with different
 values of this morphological characteristics can consume the
 same resources and are in competition for these resources.
 It can be for example the size of the beak for birds eating the same
 seeds or the height of the herbivores eating the same plants.


We will show that the homogeneous equilibrium stable in the case of the 
usual logistic equation
can become unstable in the case of nonlocal consumption of resources.
It is a new mechanism of pattern formation where spatial structures can 
appear not because of the
competition between two species as in the case of Turing structures but 
for a single species.
In the context of population dynamics it can be intepreted as emergence 
of new species due to
the intra-specific competition and random mutations.
If the population is initially localized it will split into several 
separate subpopulations.
Its dynamics resembles Darwin's schematic representation of the evolution 
process.
Such behavior cannot be described by the logistic equation without nonlocal 
terms.

Similar to the conventional reaction-diffusion equation, reaction-diffusion 
equation with integral terms
can have travelling wave solutions.
Depending on the parameters it can be either a usual travelling wave with a 
constant speed
and a stable homogeneous equilibrium behind the wave or a periodic wave with a 
periodic in space
solution behind the wave.
Their speeds of propagation (average speed for the periodic wave)
are the same as the minimal speed for the monostable reaction-diffusion 
equation.
It is determined only by the diffusion coefficients and by the derivative 
of the nonlinearity at the
unstable equilibrium. 

We note that there is a number of works where the diffusion term is replaced
by an integral term describing a nonlocal interaction (see, e.g., [4]).
In spite of some similarity of the models they are quite different.
The scalar equation with the nonlocal diffusion term satisfies
the comparison principle and cannot describe emergence of spatial
patterns and propagation of periodic waves.

Mathematical properties of some integro-differential evolution
equations are studied in [5,6].



 \vspace*{0.5cm}

 \setcounter{section}{2}
 \setcounter{equation}{0}

 \centerline{\bf 2. Model with a nonlocal interaction}

 \vspace*{0.5cm}

 We begin with a well established in population dynamics model
 taking into account migration of individuals.
 We suppose that they move randomly during their life time
 and all direction of motion are equally possible.
 If $p(r) dr$ is the probability of the displacement on
 the distance between $r$ and $r + dr$ at the unit time measured
 in generations, then
 $$ \rho = \sqrt{\int_0^\infty r^2 p(r) dr} $$
 is a meansquare displacement during one generation (see, e.g.
 [7]).
 Assuming that migration is independent of the rate of birth and
 death, we obtain a local change of the density $N$ of the
 population at the point $(x,y)$ during the time interval $\Delta
 t$:
 $$ \Delta N(x,y,t) = \left( \int N(x',y',t)
 \frac{p(r)}{2 \pi r} dx' dy' - N(x,y,t) + F(N,x,y) \right)
 \Delta t , $$
 where $r = \sqrt{(x'-x)^2+(y'-y)^2}$.
 Using the Taylor expansion for $N(x',y',t)$ around $(x,y)$ we
 obtain the nonlinear diffusion equation
 \begin{equation}
 \label{2.1}
  \frac{\partial N}{\partial t} = d  \left(
 \frac{\partial^2 N}{\partial x^2} +
 \frac{\partial^2 N}{\partial y^2} \right) + F ,
 \end{equation}
 where $d = \rho^2/4$ is the diffusion coefficient,
 $F$ is the local rate of growth of the population.
 Usually it is considered in the form
 \begin{equation}
 \label{2.1a}
 F = (B - D) N ,
 \end{equation}
 where $B$ is the birth function and $D$ is the death function.
 The birth function $B(N)$ has a sigmoidal form and it is constant
 for $N$ sufficiently large.
 The death function $D$ is usually considered in the form
 \begin{equation}
 \label{2.2}
 D(N) = b + k N ,
 \end{equation}
 where the first term in the right-hand side describe the natural
 death, and the second term its increase because of the
 competition for limited resources.
 This form of the death function implies that the resources
 are consumed locally, that is exactly at the same point $(x,y)$ where
 the individual is located.

 In a more general and more realistic case the resources can
 be consumed in some neighborhood of this point.
 We take this into account in the death function:
 $$ D(x,y,t) = b + k \int_0^t
 \left(\int_{(x,y)} \phi(x-x',y-y',t-t') N(x',y',t') dx' dy'\right) dt' . $$
 Here $\phi(x-x',y-y',t-t')$ is the function that shows how
 the individuals located at $(x',y')$ and at time $t'$
 influence the resources at the point $(x,y)$ at time $t$.
 This expression for the death function means that survival
 of individuals at the point $(x,y)$ depends on the resources
 available at this point.
 On the other hand, these resources are consumed not only by the
 individuals that have their average position at this point but
 also by the individuals from some area around it.

 If the motion of the individuals is fast with respect to the
 given time scale, then the location of the individual
 with an average position at the origin is time independent
 and given by the probability density function $\phi_0(x,y)$.
 If we assume moreover that the resources are renewed
 with a small characteristic time compared with other time
 scales (reproduction, death, migration), then the influence
 function $\phi(x-x',y-y',t-t')$ is concentrated at $t=t'$
 and equals zero for $t' < t$:
 $$ \phi(x-x',y-y',t-t') = \phi_0(x-x',y-y) \delta(t-t') . $$
 Thus we obtain the modified expression of the death function:
 \begin{equation}
 \label{2.3}
 D(x,y,t) = b + k \int_{(x,y)} \phi_0(x-x',y-y') N(x',y',t) dx' dy' .
 \end{equation}
 We note that if $\phi_0(x,y)$ is
 the Dirac $\delta$-function, then it is reduced to (\ref{2.2}).
 Substituting (\ref{2.1a}) and (\ref{2.3})  into (\ref{2.1})
 we obtain the integro-differential equation
 \begin{equation}
 \label{1.0}
 \frac{\partial N}{\partial t} = d  \left(
 \frac{\partial^2 N}{\partial x^2} +
 \frac{\partial^2 N}{\partial y^2} \right) +
 \left( B(N) -  b - k \int \phi_0(x-x',y-y') N(x',y',t) dx' dy'
 \right) N .
 \end{equation}
 We will analyze this model in the one-dimensional
 spatial case where $B(n)$ is a constant.

\vspace*{1cm}



 \setcounter{section}{3}
 \setcounter{equation}{0}

 \centerline{\bf 3. Evolution of species}

 \vspace*{0.5cm}

 In the previous section we considered a spatial distribution
 of individuals in a population and its evolution in time.
 Let us characterize the population not by individuals but by
 genes.
 Consider a gene $G$ and suppose that due to mutations it can
 be in different states, which we denote by $x_i$, $i=1,...,N$.
 Let $c_i(t)$ be the number of genes in the population at the
 state $x_i$ at time $t$.
 We suppose further that mutations can happen only between the
 neighboring states, that is for the state $x_i$ the gene can
 go to one of the states $x_{i-1}$ and $x_{i+1}$.
 If we assume that all transitions are reversible and have equal
 probability $p$, then the number of genes changing their state
 from $x_{i-1}$ to $x_i$ during the unit time will be $p c_{i-1}$.
 The number of opposite transitions is $p c_i$.
 Therefore the flux of the concentration at $x_i$ is
 $$ f_i = p (c_{i-1} - c_i) . $$
 The probability $p$ of the transition depends on the distance
 between the states and will be considered in the form
 $$ p = \frac{d}{x_i - x_{i-1}} , $$
 where $d$ is some constant.
 We suppose for simplicity that the distances between the
 neighboring states are the same and that $d$ does nor depend
 on $i$.
 The last expression is well defined if the distance between the
 states is sufficiently large.

 Thus, under the assumptions above we come to the usual diffusion
 fluxes.
 Passing from the discrete to continuous model we obtain the
 equation
 \begin{equation}
 \label{6.1}
 \frac{\partial c}{\partial t} = d  \frac{\partial^2 c}{\partial x^2}
 + (B - D) c ,
 \end{equation}
 where $c(x,t)$ is the concentration of genes at the state $x$ at
 time $t$, and, as above, $B$ is a birth function,
 $D$ is a death function.
 Taking into account the nonlocal consumption of the resources we
 come to the same model as in Section 2.

 \vspace*{1cm}


 \setcounter{section}{4}
 \setcounter{equation}{0}

 \centerline{\bf 4. Stationary solutions}

 \vspace*{0.5cm}

 Consider the equation
 \begin{equation}
 \label{1.1}
 \frac{\partial c}{\partial t} = d  \frac{\partial^2 c}{\partial x^2}
 + c \left( a - \int_{-\infty}^\infty \phi(x-y) c(y) dy \right) - b c
 \end{equation}
 for $x \in \mathbb R$.
 Here $a$ and $b$ are some constants, $\phi(y)$ is a function
 with a bounded support and such that
 \begin{equation}
 \label{1.3}
 \int_{-\infty}^\infty \phi(y) d y = 1 .
 \end{equation}
 We will assume that $\phi(y)$ is nonnegative, even and identically
 equal zero for $|y| \geq N$.
 If its support tends to zero, then we obtain in the limit the
 $\delta$-function, and equation (\ref{1.1}) becomes the usual
 reaction-diffusion equation
 \begin{equation}
 \label{1.2}
 \frac{\partial c}{\partial t} = d  \frac{\partial^2 c}{\partial x^2}
 + c ( a - c) - b c
 \end{equation}
 often considered in population dynamics and in other
 applications.

 It can be easily verified that equation (\ref{1.1}) has
 two homogeneous in space stationary solutions, $c_0(x)=0$ and
 $c_1(x)=a-b$.
 These are also stationary solutions of equation
 (\ref{1.2}).

\vspace*{1cm}


 \setcounter{section}{5}
 \setcounter{equation}{0}

 \centerline{\bf 5. Stability}

 \vspace*{0.5cm}

 If $a > b$, then $c_1(x)$ is a stable stationary solution of
 equation (\ref{1.2}).
 We will show that it can be unstable as a solution of
 equation (\ref{1.1}).
 Linearizing this equation about $c_1$, we obtain the eigenvalue
 problem
 \begin{equation}
 \label{4.1}
  d  c^{''} - \sigma  \int_{-\infty}^\infty \phi(x-y) c(y) dy =
  \lambda c ,
 \end{equation}
 where $\sigma = a-b$.
 Applying the Fourier transform to this equation, we have
 \begin{equation}
 \label{4.2}
 ( d \xi^2 + \sigma \tilde \phi(\xi) + \lambda) \tilde c(\xi) = 0,
 \end{equation}
 where the tilde denotes the Fourier transform of the
 corresponding function.
 This equation has nonzero solutions at the points of the
 spectrum, that is for such values of $\lambda$ that
 \begin{equation}
 \label{4.3}
 d \xi^2 + \sigma \tilde \phi(\xi) + \lambda = 0 .
 \end{equation}
 Consider the function
 $$ \Phi(\xi) = d \xi^2 + \sigma \tilde \phi(\xi) . $$
 If it becomes negative for some values of $\xi$, then
 there are positive $\lambda$ satisfying equation (\ref{2.3}).
 In this case, the stationary solution $c_1$ is unstable.

 We note that
 $$ \tilde \phi(0) = \int_{-\infty}^\infty \phi(y) dy = 1 . $$
 Therefore, for any given continuous  and bounded function
 $\tilde \phi(\xi)$, $\Phi(\xi)$ becomes strictly positive
 if $d$ is sufficiently large.
 Thus, there is no instability if the diffusion coefficient
 is large.

 By virtue of the assumption that the function $\phi(y)$ is even,
 we have
 $$ \tilde \phi(\xi) = \int_{-\infty}^\infty \phi(y) \cos(\xi y)
 dy. $$
 The assumption that the function $\phi(y)$ has a bounded support
 and condition (\ref{1.3})
 are not essential for the stability analysis.
 For the functions
 $$ \phi_1(y) = e^{-ay^2} , \;\; \phi_2(y) = e^{-a|y|} \;\; (a>0),$$
 we have
 $$ \tilde \phi_1(\xi) = \sqrt{\frac{\pi}{a}} \; e^{-\xi^2/(4a)}
 , \;\;
 \tilde \phi_2(\xi) = \frac{2a}{a^2 + b^2} , $$
 respectively.
 Therefore, $\Phi(\xi)$ is positive for all $\xi$, and the
 stationary solution $c_1$ is stable.

\vspace*{1cm}


 \setcounter{section}{6}
 \setcounter{equation}{0}

 \centerline{\bf 6. Piece-wise constant distribution}

 \vspace*{0.5cm}

 Consider the following function:
\begin{equation}
 \label{5.1}
\phi(y) = \left\{
 \begin{array}{ccc}
 1/(2 N) & , &  |y| \leq N \\
 0       & , &  |y| > N
 \end{array} .
 \right.
\end{equation}
 Then
 $$ \tilde \phi(\xi) = \frac{1}{2N} \; \int_{-N}^N \cos(\xi y) dy =
 \frac{1}{\xi N} \;  \sin(\xi N)  $$
 and
 $$ \Phi(\xi) = d \xi^2 +  \frac{\sigma}{\xi N} \;  \sin(\xi N) . $$
 Depending on the parameters, this function can be strictly
 positive or can change sign.
 It can be easily verified that $\Phi(0) = \sigma > 0$, and
 it is positive for $|\xi|$ sufficiently large.
 If the diffusion coefficient is large, then $\Phi(\xi)$ is
 positive for all $\xi$, and the stationary solution is stable.
 The instability can appear if $d$ is small enough.

 In the critical case, in which the function $\Phi(x)$ is
 nonnegative and not strictly
 positive, it satisfies the equalities
 $$ \Phi(\xi_0) = 0 , \;\;  \Phi'(\xi_0) = 0 $$
 for some $\xi_0$.
 It follows from these relations:
 $$ d \xi_0^2 +  \frac{\sigma}{\xi_0 N} \;  \sin(\xi_0 N) = 0 , \;\;
 2 d \xi_0 - \frac{\sigma}{\xi_0^2 N} \;  \sin(\xi_0 N) +
 \frac{\sigma}{\xi_0} \;  \cos(\xi_0 N) = 0 . $$
 Let us introduce the following notations:
 $\mu = d/\sigma, z = \xi_0 N$.
 Then the previous equations can be written in the form
 \begin{equation}
 \label{3.1}
 \mu \xi_0^2 + \frac{\sin z}{z} = 0, \;\;\;
 2 \mu \xi_0^2 - \frac{\sin z}{z} + \cos z = 0 .
 \end{equation}
 Eliminating $\mu \xi_0^2$, we obtain from these two equations
 the equation with respect to $z$:
 \begin{equation}
 \label{3.2}
 \tan z = \frac13 z .
 \end{equation}
 This equation has an infinite number of solutions.
 We will consider only positive solutions.
 Denote them by $z_1, z_2,...$
 For each of them we can find from the first equation in
 (\ref{3.1}) the relation
 \begin{equation}
 \label{3.3}
 \mu = - N \frac{\sin z_j}{z_j^3} , \;\; j=1,2,...
 \end{equation}
 Thus we have  a countable number of curves on the
 $(N,\mu)$-plane.
 Each of them gives the values of parameters for which the
 spectrum of the eigenvalue problem (\ref{2.1}) crosses
 the origin.
 We note that $\sin z_j$ is positive for even $j$ and negative
 for odd $j$.
 Since $\mu$ is positive, we will consider only odd values of $j$.

 It can be easily verified that
 $$ - \frac{\sin z_1}{z_1^3} > - \frac{\sin z_3}{z_3^3} > ... $$
 Therefore
 $$ \mu_1(N) > \mu_3(N) > ... \;\;\;\; {\rm for} \;\;\; N > 0 . $$

 We have seen in the previous section that there is no instability
 if $d$ is sufficiently large.
 Therefore, $\mu > \mu_1(N)$ corresponds to the stability region,
 $\mu < \mu_1(N)$ corresponds to the instability region.
 The stability boundary is given by the equation
 \begin{equation}
 \label{3.4}
 \mu = \mu_1(N) .
 \end{equation}
 We can also find the
 period $\tau$ of the spatial structure at the stability boundary:
 $$ \tau = \frac{2\pi}{\xi_0} = \frac{2 \pi N}{z_1} . $$

 Let us consider a fixed value of $\mu$ and vary $N$.
 For $N$ small enough the homogeneous stationary solution is
 stable.
 For some critical value of $N$ such that (\ref{3.4}) is
 satisfied, it loses its stability.
 For greater values of $N$, the solution remains unstable
 and other points of the spectrum can cross the origin.

 Suppose now that we consider a bounded interval with the length
 $L$ and with the periodic boundary conditions.
 The analysis above remains valid but we need to take into account
 the period of the perturbation.
 We fix $\xi = 2 \pi/L$ and consider the equation $\Phi(\xi)=0$.
 We obtain the equation
 $ \sin z = - \mu \xi^2 z . $
 If $\mu$ is sufficiently large, then the equation does not have
 nonzero solutions.
 In the critical case we should also equate the derivatives of
 the function in the left-hand side and in the right-hand side of
 the last equality:
 $ \cos z = - \mu \xi^2 . $
 Therefore we can find $z$ from the equation
 $ \tan z = z , $
 and the critical value of $\mu$ from the relation
 $ \mu = - (\sin z)/(\xi^2 z) . $

\vspace*{1cm}


 \setcounter{section}{7}
 \setcounter{equation}{0}

 \centerline{\bf 7. Emergence of structures}

 \vspace*{0.5cm}

Equation (\ref{2.1}) is considered on the segment $[0,L]$ with
periodic boundary conditions. It is discretised by using an
explicit finite differences scheme. The function $\Phi$ is
piece-wise constant as described in equation (\ref{5.1}). The
integral is approximated by the trapeze method.


\begin{figure}
\epsfxsize=12cm
\epsfysize=6cm
\epsfbox{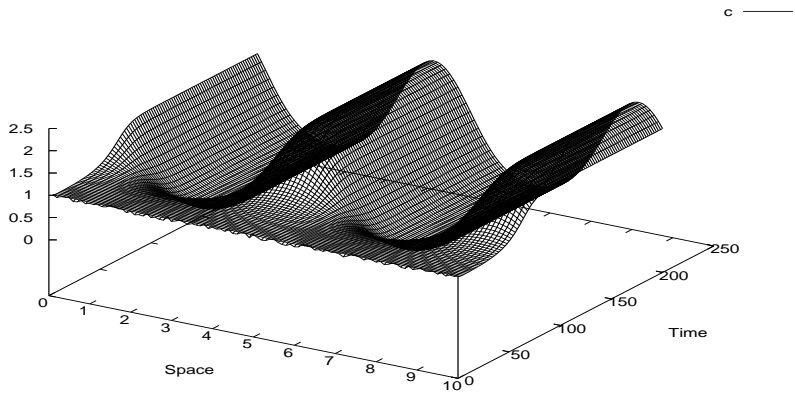}
\caption{Emergence of spatial structure}
\label{fig1}
\end{figure}

 Figure 1 
 shows the behavior of the solution of equation (\ref{1.1})
 for the values of parameters
 $d=0.05$, $a=2$, $N=3$, $b=1$ and $L=10$
 and for the initial condition
 $c(x,0)=a-b$ perturbed by a small random noise.
 This initial homogeneous equilibrium is unstable and a spatial structure
 appears.
 Note that in accordance with the stability analysis (Section 6)
 the homogeneous equilibrium is
 unstable only for small values of the diffusion coefficient $d$.
 For the same values of the parameters above and  $d>0.11$, the homogeneous
 equilibrium is stable.

\begin{figure}
\epsfxsize=12cm
\epsfysize=6cm
\epsfbox{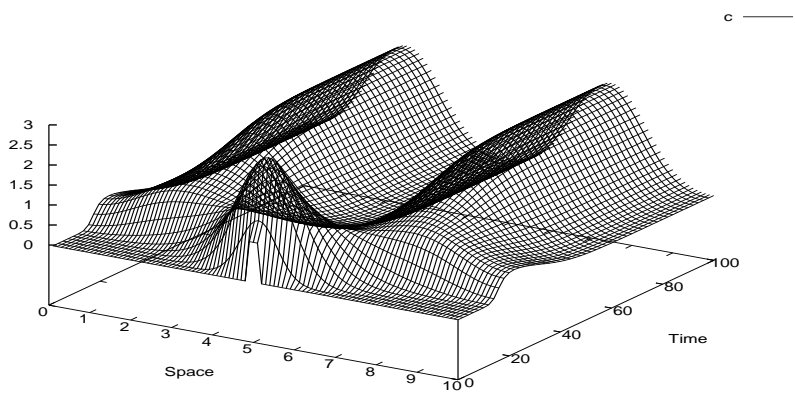}
\caption{Splitting of a population}
\label{fig2}
\end{figure}

 We consider next the same values of parameters as for Figure 1
 and the initial condition $c(0,x)=a-b$ in a small interval centered
 at the middle of $[0,L]$ and $c(0,x)=0$ otherwise.
 At the first stage of its evolution the solution grows and
 diffuse slightly (Figure 2),
 its maximum remains at the middle of interval $[0,L]$.
 At the second stage the solution decreases at the center and
 grows near the borders of the interval.
 We observe here that the intra-specific competition results in
 appearance of two separated sub-populations symmetrically situated around the middle of
 $[0,L]$.

\begin{figure}
\epsfxsize=18cm
\epsfysize=9cm
\epsfbox{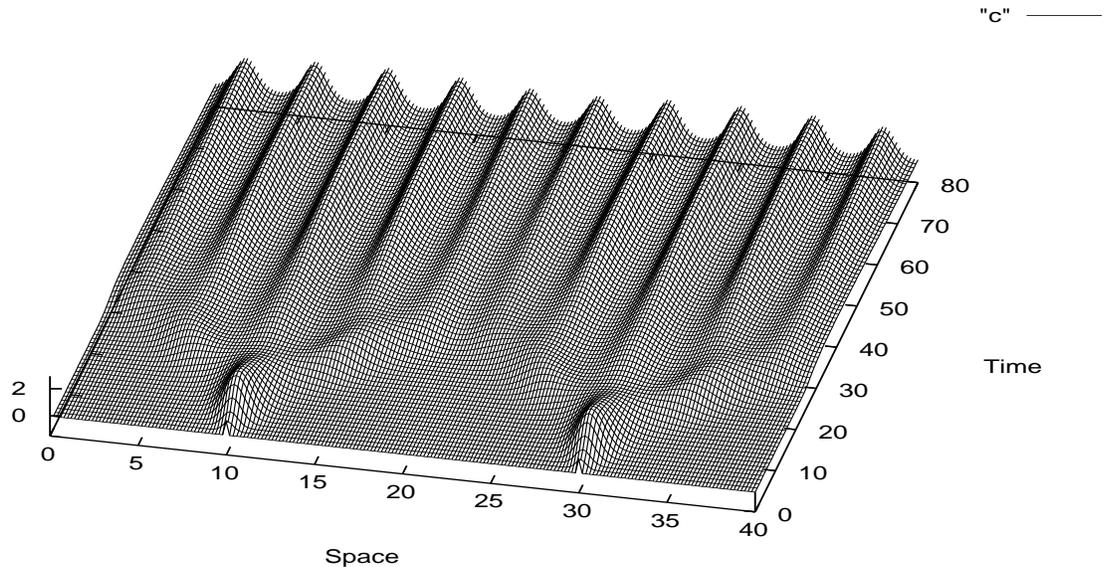}
\caption{Evolution}
\label{fig3}
\end{figure}

 The case of a larger interval with two initially separated
 sub-populations is shown in Figure 3.
 Each of them first splits into three sub-populations.
 Those in the center do not change any more.
 The sub-populations from the sides split once more into
 two sub-populations each of them.
 These consecutive splitting into sub-populations has an
 interpretation related to the evolution of species.
 We discuss this question below.

 The numerical simulations shown in Figure 4 are carried out
 for the same values of the parameters except for the diffusion
 coefficient which is smaller in this case.
 The consecutive sub-populations are almost disconnected here.
 Each new structure emerges from very small values of the
 concentration.

 We consider finally the case where the function $\phi(y)$ is
 not symmetric.
 It is similar to that given by (\ref{5.1}) but equals $0$
 for $y < 0$ (see Section 9).
 The dynamics of the solution is different in this case:
 the first maximum of the concentration moves to the left, then
 the second maximum appears and moves remaining at fixed
 distance from the first one, and so on (Figure 5).

\begin{figure}
\epsfxsize=18cm
\epsfysize=9cm
\epsfbox{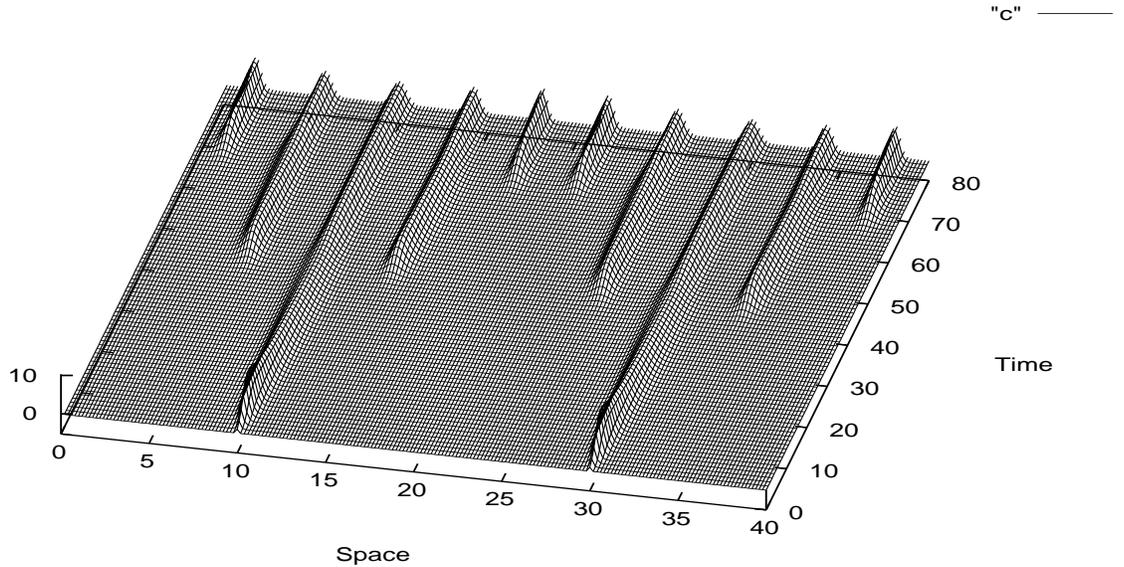}
\caption{Small diffusion coefficient}
\label{fig4}
\end{figure}

\begin{figure}
\epsfxsize=18cm
\epsfysize=9cm
\epsfbox{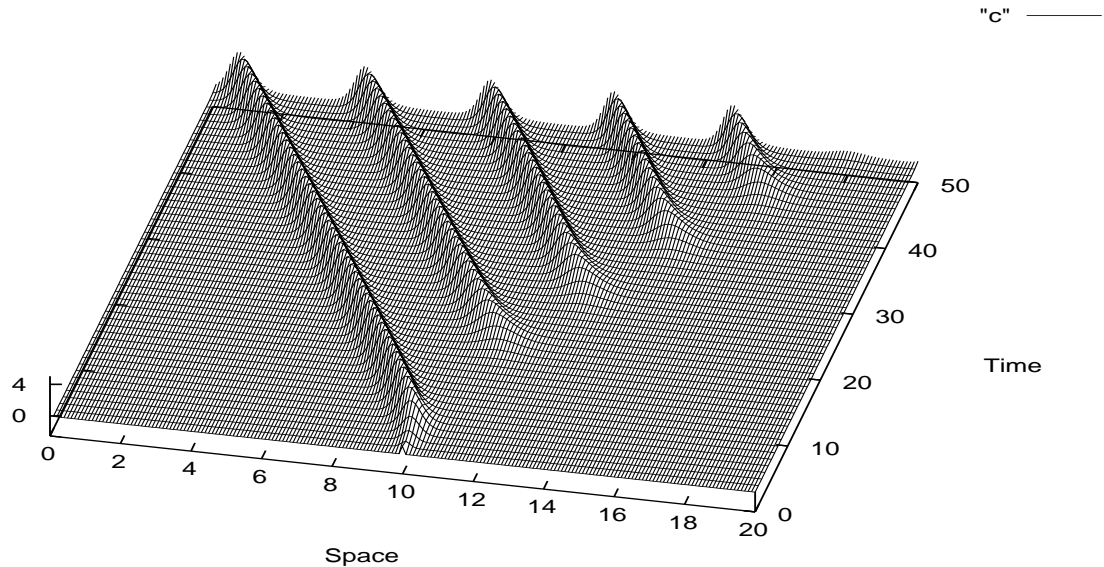}
\caption{Asymmetric evolution}
\label{fig5}
\end{figure}


\vspace*{1cm}

 \setcounter{section}{8}
 \setcounter{equation}{0}

 \centerline{\bf 8. Approximate solution}

 \vspace*{0.5cm}

 In this section we consider an approximation of equation
 (\ref{1.1}), which is not well justified mathematically but
 allows us to explain some of its properties.
 Consider the integral in the right-hand side of this equation
 $$ \int_{-\infty}^\infty \phi(x-y) c(y) dy =
 \int_{-\infty}^\infty \phi(y) c(x-y) dy $$
 (the time dependence of the function $c$ is for brevity omitted).
 We will use the Taylor expansion
 $$ c(x-y) = c(x) - c^{'}(x) y + \frac12 c^{''}(x) y^2 + ... $$
 If $\phi$ is an even function then we obtain
 \begin{equation}
 \label{7.1}
 \int_{-\infty}^\infty \phi(x-y) c(y) dy =
 c(x) + \gamma c^{''}(x) + ... ,
 \end{equation}
 where
 $$ \gamma = \frac12 \int_{-\infty}^\infty \phi(y) y^2 dy . $$
 Taking only the first term in the right-hand side of expansion
 (\ref{7.1}), we obtain the reaction-diffusion equation
 (\ref{1.2}) instead of the integro-differential equation
 (\ref{1.1}).
 Equation (\ref{1.2}) is well-studied.
 It is known that it has a travelling wave solution, that is
 the solution of the form $c(x,t) = u(x-st)$, where $s$
 is the wave velocity.
 In fact, monotone in $x$ waves for this equation exist for all
 $s \geq s_0 = 2 \sqrt{d \sigma}$, where $\sigma = a-b$.
 If the nonlinearity $F(c)$, which in this particular case
 equals $c(\sigma-c)$, satisfies the condition $F'(c) \leq F'(0)$,
 then the minimal speed $s$ is determined by $F'(0)$ and
 not by the nonlinearity itself.
 In this sense the minimal speed is stable with respect to
 the perturbation of the nonlinearity.
 It is also known that the solutions of the parabolic equation
 with localized in space initial conditions converge to the
 wave with the minimal velocity [8].

 Consider next the two-term approximation in (\ref{7.1}).
 We obtain the equation
 \begin{equation}
 \label{7.2}
 \frac{\partial c}{\partial t} = d  (1- \gamma c) \frac{\partial^2 c}{\partial x^2}
 + c (\sigma - c) .
 \end{equation}
 The corresponding stationary equation
 $$ d  (1- \gamma c) c^{''} + c (\sigma - c) = 0 $$
 has periodic in $x$ solutions if $b \sigma > 1$.
 Indeed, it can be written as the system of two first order
 equations:
 $$ c' = p, \;\; p' = - \Phi(c) , $$
 where $\Phi(c) = c (\sigma - c)/( d  (1- \gamma c))$.
 Since $\Phi(\sigma)=0$, $\Phi'(\sigma) > 0$, then
 there exists a family of limit cycles around the stationary
 point $c=\sigma,p=0$ of this system.

 Thus, if $\gamma=0$ there are travelling waves with the limits
 $c=0$ and $c=\sigma$.
 Each of these values is a stationary solution of the evolution
 equation, the first one being unstable while the second one
 is stable.
 The minimal speed is determined by the linearized about
 $c=0$ equation.
 If $\gamma \neq 0$, the same homogeneous stationary
 solutions exist but $c=\sigma$ can become unstable.
 Therefore the wave can connect the homogeneous solution
 $c=0$ at one infinity with the periodic solution at another
 infinity.
 However, it cannot be done for equation (\ref{7.2}) because
 the coefficient of the second derivative will change sign.
 We need to consider the complete equation (\ref{1.1}).
 Its linearizations about $c=0$ and about $c=\sigma$
 have the same properties as for equation (\ref{7.2}).
 Hence we can expect the existence of a periodic in time travelling wave,
 that is
 $$ c(x,t) = u(x-st,t), $$
 where $u(x,t)$ is a periodic in $t$ solution of the corresponding
 equation.
 If we assume that the minimal speed is determined by the
 linearized about $c=0$ equation, as before, then it will
 have the same value $s_0 = 2 \sqrt{d \sigma}$.
 For the values $d = 0.06, a = 2, b = 1$ of the parameters
 considered in the simulations on Figure 6, we obtain
 $s_0 \approx 0.48$.
 This result is in a good agreement with the speed of propagation
 found numerically.

\begin{figure}
\epsfxsize=16cm
\epsfysize=8cm
\epsfbox{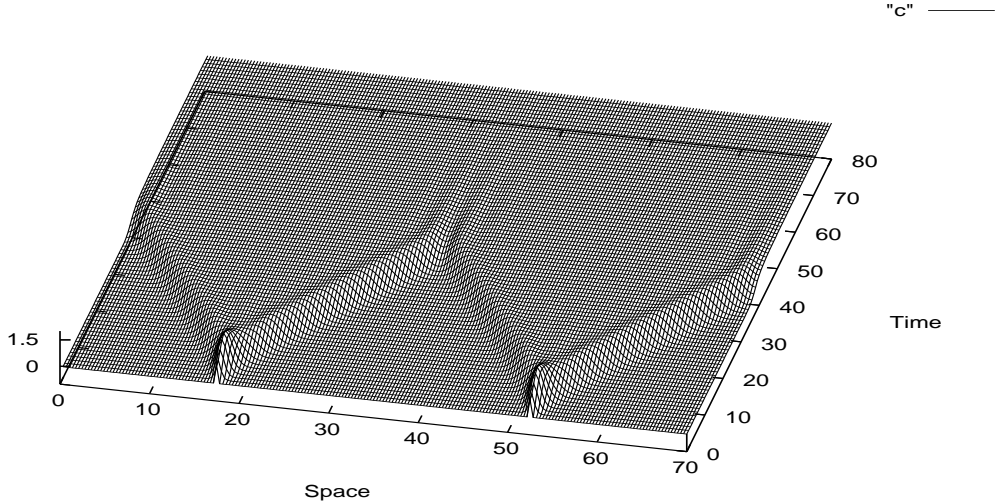}
\caption{Stable equilibrium behind the wave}
\label{fig6}
\end{figure}

\begin{figure}
\epsfxsize=16cm
\epsfysize=8cm
\epsfbox{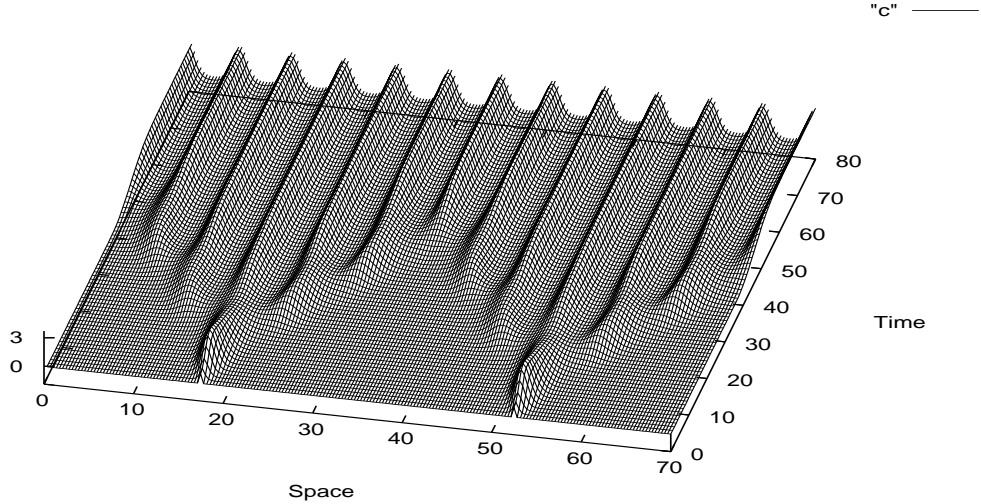}
\caption{Periodic wave propagation}
\label{fig7}
\end{figure}

It is interesting to compare the results on Figures 6 and 7.
The values of the parameters are the same except for $N$ which
determines the width of the support of the function $\phi$.
In the second case (Figure 7) it is twice more than in the first
case (Figure 6).
In agreement with the linear stability analysis the equilibrium
behind the wave is stable in the first case and unstable
in the second one.
However the speeds of wave propagation in both cases
is the same.
This confirms the assumption that the minimal speed is determined
by the system linearized about $c=0$, and that the solution
converges to the wave with the minimal speed.

\vspace*{1cm}



 \setcounter{section}{9}
 \setcounter{equation}{0}

 \centerline{\bf 9. Asymmetric evolution}

 \vspace*{0.5cm}

 In Section 6 we considered the case where the function $\phi$
 was even.
 Therefore its Fourier transform $\tilde \phi$ was a real-valued
 function, and the principal eigenvalue crossed the imaginary
 axis through the origin.
 On the other hand, the numerical simulations presented in Section
 7 showed that behavior of solutions in the case of the asymmetric
 function $\phi$ could be different in comparison with the symmetric
 case.
 Consider for certainty the following function:
 \begin{equation}
 \label{9.1}
 \phi(y) = \left\{
 \begin{array}{ccc}
 1/ N & , &  0 \leq y \leq N \\
 0       & , &  \mbox{otherwise}
 \end{array} .
 \right.
 \end{equation}
 Then
 $$ \tilde \phi(\xi) = \frac{1}{N} \;
 \int_{0}^N (\cos(\xi y) + i \sin(\xi y))dy =
 \frac{1}{\xi N} \;  \sin(\xi N)
 +  \frac{i}{\xi N} \; ( \cos(\xi N)- 1) . $$
 Hence the real part of the function $\Phi(\xi)$ is the same as
 before but the imaginary part is now different from zero.
 This means that the stability boundary is the same as for the
 symmetric function but the behavior of solutions is different.

  To describe the behavior of solutions at the stability
  boundary consider the linearized equation
 \begin{equation}
 \label{9.2}
 \frac{\partial c}{\partial t} = d \;
 \frac{\partial^2 c}{\partial x^2} \;
 - \sigma  \int_{-\infty}^\infty \phi(x-y) c(y) dy .
 \end{equation}
 We look for its solution in the form
 $$ c(x,t) = \cos (\xi (x-st)) . $$
 Substituting it in (\ref{9.2}) and separating the real and
 imaginary parts, we obtain
 $$ d \xi^2 + \sigma \int_{-\infty}^\infty \phi(y) \cos (\xi y) dy = 0,
 $$
 $$ \xi s = - \int_{-\infty}^\infty \phi(y) \sin (\xi y) dy . $$
 The first equality is satisfied at the stability boundary.
 The second equality gives the relation between the speed
 $s$ and the frequency $\xi$.
 For the function (\ref{9.1}) we have
\begin{equation}
\label{9.3}
  s = \frac{1}{\xi^2 N} \; (\cos (\xi N) - 1) . 
\end{equation}

\begin{figure}
\epsfxsize=18cm
\epsfysize=9cm
\epsfbox{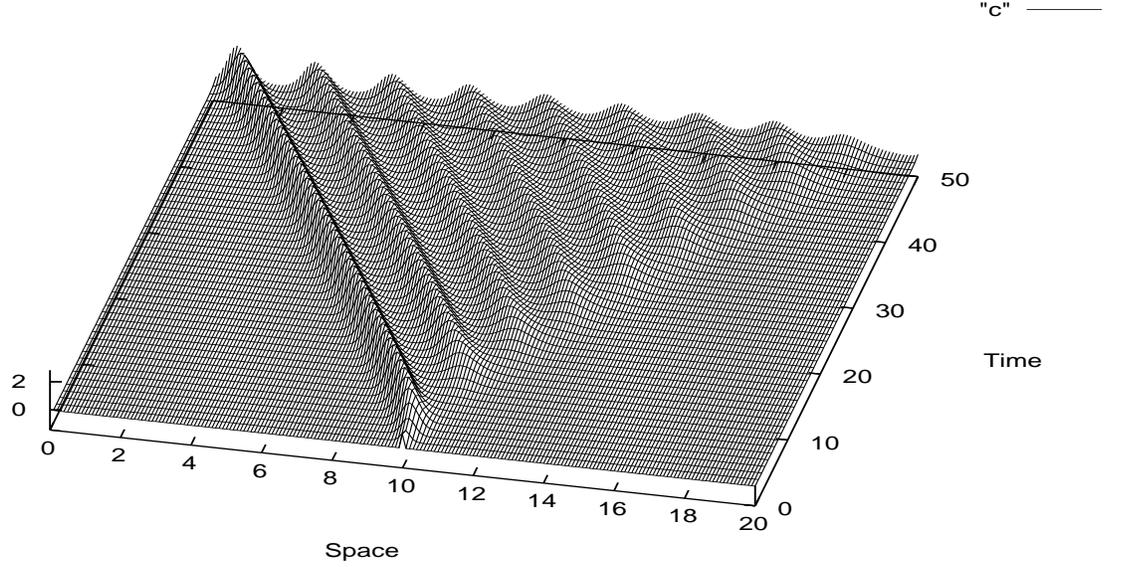}
\caption{Asymmetric evolution}
\label{fig8}
\end{figure}

We verify this relation with the numerical simulations shown
in Figure 8.
The values of parameters here ($d=0.01, \sigma=1, N=1.1$)
are such that the homogeneous solution is unstable but they
are not far from the stability boundary ($N=1$).
The speed $s$ of propagation equals the minimal speed
$s_0 = 2 \sqrt{d \sigma} = 0.2$.
The value of $\xi$ taken from the numerical results is approximately
$2.85$.
The right-hand side in (\ref{9.3}) equals approximately $0.22$.
Therefore this equality is approximately satisfied.
We note that it is satisfied at the stability boundary but
not outside it. 
For larger values of $N$ the frequency of oscillations
is quite different (see Figure 5, where all the parameters are the same
as in Figure 8, except that $N=2$) though the speed of propagation
remains the same.

\vspace*{1cm}



 \setcounter{section}{10}
 \setcounter{equation}{0}

 \centerline{\bf 10. Discussion}

 \vspace*{0.5cm}

 In this work we study population dynamics with nonlocal consumptions
 of resources.
 We consider the conventional logistic equation with diffusion and with
 a nonlinear term describing reproduction and mortality of individuals
 in a population.
 The difference in comparison with the classical in the population dynamics model is
 that the death function depends not only on the density of the individuals
 at the given spatial point but also on their density in some its neighborhood.
 This is related to the assumption that the individuals in the population consume
 the resources not only at the point where they are located but also in some
 area around this point.

 The nonlocal consumption of the resources results in appearance of spatial
 structures.
 The homogeneous equilibrium, which would be stable otherwise, becomes
 unstable.
 This is a new mechanism of pattern formations different with respect to
 Turing structures.
 As it is well known, Turing or dissipative structures are described by
 reaction-diffusion systems of at least two equations for two different concentrations.
 They are often associated with an activator and inhibitor in some chemical systems.
 This mechanism can be briefly described as a short range activation and a long range
 inhibition.
 Though this mechanism is often applied to describe biological pattern formation,
 and in some cases it gives patterns very close to those observed in reality,
 it remains basically phenomenological (see [9]).

 The model considered in this work can show the emergence of patterns in a single
 population, without competition between different species.
 To our knowledge, it is a new effect for modelling in population
 dynamics.
 Its biological interpretation is that an intra-specific
 competition can result in a spatial localization of the
 individuals in the populations.
 We can suppose that formation of herds of herbivorous animals
 or formation of local communities in human populations, though
 much more complex, can be related to this mechanism.

 The possibility of the appearance of spatial structures depends
 on the function $\phi(y)$.
 In particular, the piece-wise constant function introduced in
 Section 6 can lead to the instability but not the error function
 (Section 5).
 Let us discuss a possible biological meaning of the function
 (\ref{5.1}).
 Consider an individual with his average location at $y=0$
 where he has his family life (reproduction).
 He consumes resources during the day and returns to the family
 during the night.
 Suppose that he needs 8 hours per day to consume the sufficient
 amount of resources and that he can spend 12 hours per day
 outside the nest.
 Then he has 2 hours to go to his work (consumption) place and
 2 hours to return back.
 If we neglect the cost of the displacement and if we assume that
 the working places are equally distributed, then he will go
 with a constant probability to any place in the two hours
 distance from the nest but not further than that.
 This gives us the piece-wise constant distribution  (\ref{5.1}).
 If the diffusion coefficient in this case is sufficiently small,
 then a nonhomogeneous spatial distribution will emerge.
 We should distinguish here a fast motion home-work, which is
 considered as instantaneous displacement, and a slow
 motion due to diffusion.
 The latter is in fact a displacement of the nest.
 For example, it can be a relatively slow motion to move the
 home closer to the work place.
 The characteristic time of this motion can be related to the
 development of the infrastructure and of other factors.
 If the diffusion coefficient is large, then the spatial
 structure does not appear and the homogeneous equilibrium remains
 stable.
 In human populations it can result in appearance of big
 equally populated areas (e.g. Chicago area in the USA).

\begin{figure}
\epsfxsize=16cm
\epsfysize=12cm
\epsfbox{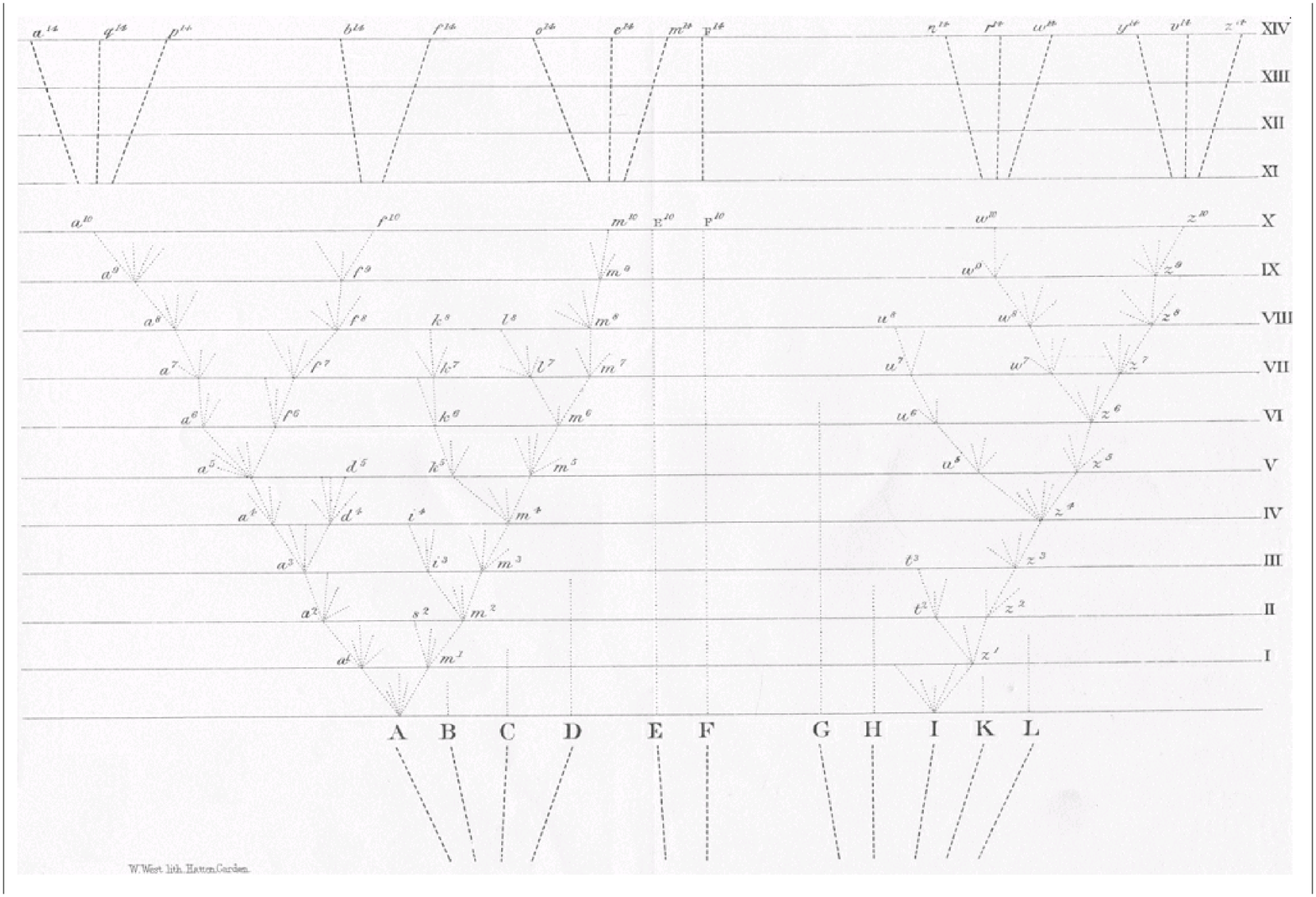}
\caption{Evolution of species according to Darwin}
\label{fig9}
\end{figure}

 On the other hand, if the cost of the displacement is high,
 then $\phi(y)$ will decrease when $|y|$ increases.
 If the influence of this factor is sufficiently strong, then
 the population will also remain equally distributed though
 for another reason: to stay close to the consumption place.

 Another context where the model considered in this work can be
 applied is evolution of species (see Section 3).
 In this case the space variable does not correspond to the
 physical space anymore but to the morphological space.
 Consider some morphological feature, for example, the size of
 bird's beak and the resources related to this feature.
 It can be grains of different sizes that birds consume.
 This brings us to the same piece-wise constant function $\phi$
 as before.
 By diffusion in this case we understand a slow change of the beak
 size from one generation to another due to random mutations
 or due to some other factors.
 The intra-specific competition for resources can result in
 the splitting of the population into one or several
 sub-populations (see Figures 2, 3) that can be considered
 as new species.
 If the diffusion coefficient is sufficiently small, then
 the new species emerge quickly in time and they are not connected
 to the previous species (Figure 4).
 This can be an explanation of the fact that transient forms in
 the evolution of species are not known.

 It is interesting to note that Figure 3 looks very close to
 the schematic representation of evolution given by Darwin
 in his ``Origin of species'' (see [10], pages 116, 117).

 Let us discuss finally the simulations shown in Figure 5.
 The function $\phi$ is not symmetric in this case.
 The biological interpretation of this assumption is that
 there is a preferential direction in the consumption of the
 resources.
 It can be related for example to a dominant direction of wind
 in the case of insects or to some other environmental or social
 factors (e.g. preferential transfer of jobs to the East and
 of workers to the West).
 In the evolution context it can be some asymmetry of
 morphological features.
 For example, if an animal can easier bent down to look for food
 than to stretch up it can result in appearance of species of a
 higher height.

\vspace*{0.5cm}


 \vspace*{0.5cm}

 \centerline{\bf References}

 \vspace*{0.5cm}

 1. G. Edelman, J. Gally.
 Degeneracy and complexity in biological systems.
 Proceedings of the National Academy of Sciences, 98 (2001), No. 24, 13763-13768.

2. J.J. Kupiec, P. Sonigo.
 Ni Dieu ni g\`ene. Pour une autre t\'eorie de l'h\'er\'edit\'e.
 Seuil, Paris, 2000.

 3. S. Atamas.
Self-organization in computer simulated selective systems.
Biosystems, 39 (1996), 143-151.

4. J. Coville, L. Dupaigne. 
 Propagation speed of travelling fronts in non local reaction-diffusion 
 equations.  Nonlinear Anal.  60  (2005),  No. 5, 797-819.

5. C. Prevost. Application des EDP
aux
probl\`emes de dynamique des populations et traitement num\'erique.
PHD, Universit\'e d'Orleans, 2004.

6. L. Desvillettes, C. Prevost, R. Ferrieres. 
Infinite Dimensional Reaction-Diffusion for Population
Dynamics.
Preprint n. 2003-04 du CMLA, ENS 
Cachan.

 7. Yu. M. Svirizhev.
 Nonlinear waves, dissipative structures, and disasters in ecology.
 Nauka, Moscow, 1987.

8. A. Volpert, Vit. Volpert, Vl. Volpert.
 Traveling wave solutions of parabolic systems.
Translation of Mathematical Monographs, Vol. 140, 
Amer. Math. Society, Providence, 1994.

 9. L. Wolpert. Principles of development. Second Edition. Oxford
University Press, Oxford, 2002.

 10. C. Darwin. 
 On the origin of species by means of natural selection. 
 John Murray, London, 1859. [1st edn].

\end{document}